\documentclass[10pt,reqno,a4paper]{article}

\usepackage{amsmath,amsthm,amsfonts,amssymb}

 \theoremstyle{definition}

\newtheorem{exmps}{Examples}
\newtheorem{alg}{Algorithm}
\newtheorem{pb}{Problem}
\newcommand{\beqa}{\begin{eqnarray*}}
\newcommand{\eeqa}{\end{eqnarray*}}

\DeclareMathOperator*{\supp}{supp}

\def\<{\left<}
\def\>{\right>}

\def\mv1{M_v^1}

\usepackage{epsf}
\usepackage{epsfig}
\usepackage{float}

\title{Nonlinear projection digital image inpainting and restoration methods}
\author{Massimo\ Fornasier }
\begin{document}
\maketitle

\begin{abstract}
This paper concerns with nonuniform sampling and interpolation methods combined with variational models for the solution of a generalized image inpainting problem and the restoration of digital signals.
In particular, we discuss the problem of reconstructing a digital signal/image from very few, sparse, and complete information and a substantial incomplete information, which will be assumed as the result of a nonlinear distortion. As a typical and inspiring example, we illustrate the concrete problem of the color restoration of a destroyed art fresco from its few known fragments and some gray picture taken prior to the damage.
Numerical implementations are included together with  several examples and numerical results to illustrate the proposed method. The numerical experience suggests furthermore that a particular system of coupled Hamilton-Jacobi equations is well-posed. 
\end{abstract}


\section{Introduction}

Imagine an important and huge art fresco, maybe dated to 1450, right at the very beginning of the Italian Renaissance, and imagine that one day a dramatic and catastrophic event destroyed it into fragments, maybe during the Second World War. It would be a very hard problem to puzzle up this fragmented opera, especially if the original dimensions were large and the number of fragments was of several thousands.\begin{figure}[ht]
\hbox to \hsize {\hfill \epsfig{file=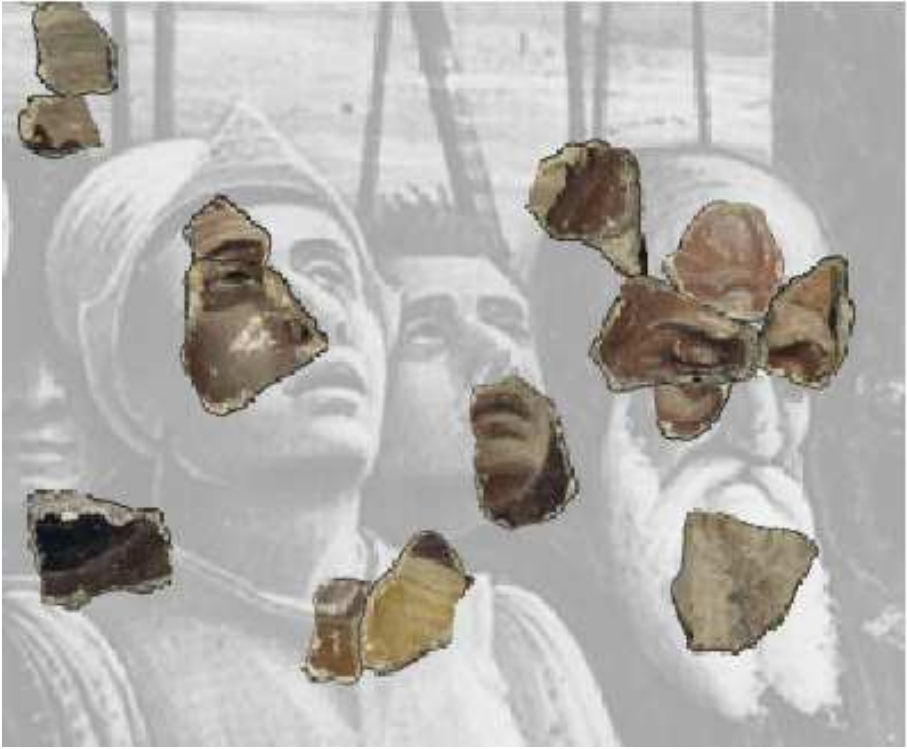,height=10cm} \hfill}
\caption{Fragmented A. Mantegna's frescoes (1452) by a bombing in the Second World War. Computer based reconstruction by using efficient pattern matching techniques}
\end{figure}

On $11^{th}$ March 1944, a group of bombs launched from an Allied airplane hit the famous Italian Eremitani's Church in Padua, destroying it together with the inestimable frescoes of Andrea Mantegna {\it et al.} contained in the Ovetari Chapel. 
People picked up and collected what was remaining of one of the most important operas of the Italian Renaissance, saving the fragments in some improvised boxes constructed with the wood of the sits of the Church. Several attempts to restore these fragments by traditional methods have been done in the last 60 years, without much success. Details on ``the state of the art'' can be found in the booklet \cite{FT2}. 
Since 1998, the author has been involved in the fascinating attempt to recall to life these frescoes, by using mathematical methods and computer based techniques. Recently a fast, robust, and efficient pattern recognition algorithm has been developed \cite{FT3,FT1} in order to detect the right position and orientation of the fragments, by means of the comparison with an old gray level image of the fresco prior to the damage\footnote{The gray images of the fresco are dated to 1920.}. This method showed to be very effective and, for example, Figure 1 illustrates some detected fragments positioned on their original place.
Unfortunately the surface covered by the original fragments is only $77$ $m^2$, while the original area was of several hundreds. This means that what we can currently reconstruct is just a {\it fraction} of what was this inestimable artwork. 
In particular, it is not known what can be the {\it original color} of those missing parts. There exist some color images of parts of the frescoes, but it is already proved that such colors are not faithful at all and the quality of such pictures is very low. So, natural questions raise: Is it possible to estimate {\it mathematically} what can be the original colors of the frescoes by using the known fragments' information? And, how {\it faithful} is this estimation? 
In this paper we want to discuss these questions and to illustrate some possible solutions.
Clearly, this problem can be reformulated in very different ways for different other situations. Therefore, the ideas developed in this context can be put into more general frames as we will illustrate in the following.
\\

The reconstruction of a {\it signal} from {\it sparse} and {\it nonuniform} sampling data is a well known problem in information theory \cite{BF,BF2}. Mathematical methods and numerical algorithms have been developed to compute missing parts of signals from few and sparse known sampling information and we refer, for example, to classical works of Feichtinger {\it et al.} \cite{FG92,FGS95,FS94} and Aldroubi/Gr\"ochenig \cite{AG01} for major details. These methods are essentially based on a {\it (quasi-)interpolation} of the signal/image by means of suitable series expansions of irregularly shifted (translated) basic functions or, in its discrete version \cite{FGS95,FS94,BG}, by series expansion of complex exponentials. 
Therefore, we may say that relevant mathematical methods and concepts maybe useful in image restoration problems are {\it interpolation} and {\it Fourier/numerical harmonic analysis} techniques. {\it Image interpolation} can be called, with a more ``artistic'' synonym, ``{\it inpainting}''. In the last years, image processing problems attracted the attention of the mathematical community, and several PDE and variational calculus techniques have been discovered to play also relevant roles to model several situations encountered in this context. We refer the reader to the nice book \cite{AK02} for a quite comprehensive description of this field. One of the first effective contributions where variational methods and PDE were used specifically for the inpainting problem is due to Beltramio, Shapiro {\it et al.} \cite{BSCB01}.
Several consequent papers came out and we refer to \cite{CS02} and \cite{AM03} as some of the most recent, interesting, and relevant for our purposes.

\begin{figure}[ht]
\hbox to \hsize {\hfill \epsfig{file=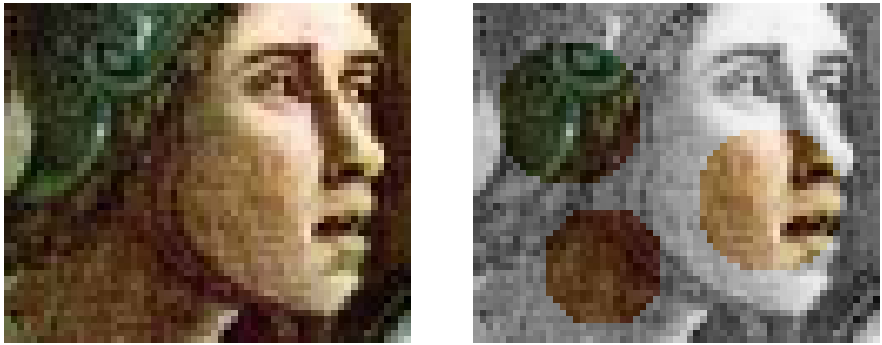,height=5cm} \hfill}
\caption{The original template picture is illustrated on the left. Few color fragments are illustrated by darker disks on the right, where  the gray level of the missing part is known. For the reader convenience, the gray level brightness has been exalted in order to distinguish it from the color part.}
\end{figure}

All these important contributions are essentially addressed to solve the problem of {\it guessing} or {\it learning} some relatively small missing part of a signal/image by using the information of the relevant known part, and they are based on deterministic methods.
The problem of learning from examples is a corner stone in information theory and artificial intelligence. Very recently a probability theory of what can be learned from a relatively large distribution of data according to a (unknown) probability measure has been formalized in the beautiful contribution by Cucker and Smale \cite{CSm}. Inspired by this work, subsequent papers, for example \cite{BCDDT}, illustrate methods to construct estimators of the probabilistically best solutions of the learning problem.
\\

The situation just exemplified by the fresco problem is slightly different.\\
We should reconstruct a {\it large} missing part, from a {\it relatively small} complete information {\it and} from a {\it partial knowledge} of the missing part itself.
In fact, see for example Figure 1, what is known of our fresco is at least the gray level of the missing part. This additional information should be exploited to recover larger parts with respect to what was possible with the classical interpolation and PDE inpainting methods only. Moreover, one should expect that the additional information can increase the probability that good estimators have to reconstruct the signal/image, in the spirit of Cucker and Smale ideas.
\\

In this paper we want to show how interpolation and deterministic PDE inpainting methods can be separately  adapted to work properly in this different situation, and how they can indeed  interplay together to achieve nice results in concrete cases. Our aim here is not to give an exhaustive description of such technique, but mainly to introduce some new framework and models. Subsequent papers will detail theoretical and numerical properties of the methods here introduced.

The paper is then organized as follows.
Section 2 describes the general framework of the problem and its particular formulation in the case of the color reconstruction of digital images from few complete samples and a large information on the gray levels. Section 3 is devoted to illustrate a simple solution to the problem, based on interpolation techniques and nonuniform sampling methods. In Section 4 we discuss variational models and evolutionary PDE to solve the problem for 1D signals and for digital images. We conclude with the formulation of a system of coupled Hamilton-Jacobi equations for the solution of the color inpainting problem, its numerical implementation, and results.

\subsection*{Acknowledgment}

The author acknowledges the financial support provided through the Intra-European Individual Marie Curie Fellowship, project FTFDORF-FP6-501018, and the hospitality of NuHAG (Numerical Harmonic Analysis Group), Department of Mathematics, University of Vienna, Austria. The author is greatly indebted with Domenico Toniolo for the work done together in the Mantegna Project and with the Mantegna Laboratory (University of Padova, Italy) that has been, with its results, the inspiration for this work.\\
A particular thank is addressed to Maurizio Falcone, Stefano Finzi Vita, Giorgio Fusco, Riccardo March, and Pierpaolo Soravia for the stimulating and nice discussions on the topic.

\subsection*{Nota on color pictures}


This paper introduces methods to recover colors in digital images. Therefore the gray level printout of the manuscript does not allow to appreciate fully the quality of the techniques illustrated and the author recommends the interested reader to access to the electronic version with color pictures available online at [...]
\\

\underline{For the review process}: Figure 15 and a corresponding animation are currently available online with colors at http://www.math.unipd.it/$\sim$mfornasi/abstractJMIV.html.

\section{Scope of the problems}

The nonuniform sampling methods to recover a band-limited function $u$ from an irregular set of sampling points $\mathcal{S}:=\{u(x_j)\}_{j \in \mathcal{J}}$ described in \cite{FG92,FGS95,FS94,BG} are based on the assumption that the sampling points are {\it dense enough} with respect to the {\it band-width} of the function to be restored. This exactly in the same spirit as the more classical and well-known Whittaker-Shannon theorem. 
This is intuitively clear in the sense that if an image is rich of details and the missing part is quite large with respect to the known part, it should be impossible to recover correctly its original version.
It is matter of the information that  one has in fact at disposal.

\begin{figure}[ht]
\hbox to \hsize {\hfill \epsfig{file=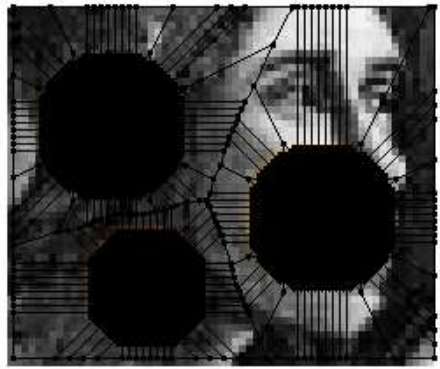,height=5cm} \hfill}
\caption{Voronoi decomposition of the domain where the nodes are the known complete color pixels.}
\end{figure}

We are interested here to discuss and to treat a relatively different {\it inverse} problem.
\begin{pb}
\label{pb1}
Consider a signal $u: \Omega \rightarrow \mathbb{R}^m$, $\Omega \subseteq \mathbb{R}^n$, and a submanifold $\mathcal{M} \subset \mathbb{R}^m$ of dimension $d \leq m$ endowed with a $C^1$ {\it nonlinear} map $\mathcal{L}:\mathbb{R}^m \rightarrow \mathcal{M}$. We call the couple $(\mathcal{M},\mathcal{L})$ a $C^1$-{\it nonlinear projection}. We assume to have a set of sampling nodes $\mathcal{N}:=\{x_j\}_{j \in \mathcal{J}} \subset \Omega$ and a subset $\mathcal{I} \subsetneq \mathcal{J}$ such that $\mathcal{C}:=\{u(x_i)\}_{i \in \mathcal{I}}$, subset of $\mathbb{R}^m$, and $\mathcal{D}:=\{\mathcal{L}(u(x_j))\}_{j \in \mathcal{J} \backslash \mathcal{I}}$, subset of $\mathcal{M}$, are known sets of values. 
The problem is then: How to reconstruct $u$ from the sampling values $\mathcal{C} \cup \mathcal{D}$.
\end{pb}

Of course, the nonlinear distortion $\mathcal{L}$ is assumed as a datum of the problem, and in concrete cases one should have estimated it a priori. This is the typical situation where {\it learning from examples} \cite{CSm,BCDDT} is necessary,  the nature of the real source of the distortion being unknown in general. In the following we will discuss briefly this estimation problem in the case of the color-gray conversion for images, exemplifying the process with the art fresco restoration (Figure 4). 

\begin{exmps}
1. Clearly if the function $\mathcal{L}_{|\text{ran}(u)}$ is injective, then it is invertible and one can construct a new set of sampling values $\mathcal{S}:= \mathcal{C} \cup \mathcal{L}_{|\text{ran}(u)}^{-1} \mathcal{D}$. If $\mathcal{S}$ is dense enough and if $u$ is a suitable band-limited function, then clearly it is possible to recover the function by the known nonuniform sampling methods.

2. If $0\in \mathcal{M}$ and $\mathcal{L}:\mathbb{R}^m \rightarrow \mathcal{M}$ has vanishing range $\text{ran}(\mathcal{L}_{|\text{ran}(u)})=\{ 0 \}$, then the problem reduces again just to the nonuniform sampling problem where $\mathcal{S}:= \mathcal{C}$.
\end{exmps}

The interesting cases are when, for example, the function $\mathcal{L}$ is not injective on $\text{ran}(u)$, $\text{dim}(\mathcal{M})<m$, and the density of the points $\{x_i\}_{i \in \mathcal{I}}$ is not enough to ensure the perfect reconstruction of the function $u$. In relevant cases, even if the information contained in $\mathcal{C}$ is not sufficient to recover the function, maybe the additional information $\mathcal{D}$ can be exploited ``somehow''.

This very general framework fits with many possible applications raising in concrete cases: For example, the recovery of a transmitted signal affected by a stationary (nonlinear) distortion, or, as in the case mentioned at the beginning of our paper, the restoration of colors of a fresco from its fragments and some gray picture taken prior to the damage.

A digital image is modeled as a function $u: \Omega\subset \mathbb{R}^2 \rightarrow \mathbb{R}^3_+$, so that, for each ``point'' $\mathbf{x}$ of the image, one associates the $3D$ vector $u(\mathbf{x})=(r(\mathbf{x}),g(\mathbf{x}),b(\mathbf{x})) \in \mathbb{R}^3_+$ of the color divided between the different channels red, green, and blue.
In particular a digitalization of the image $u$ corresponds to its sampling on a regular lattice $\tau \mathbb{Z}^2$. Let us again write $u:\mathcal{N} \rightarrow \mathbb{R}^3_+$, $u(\mathbf{x})=(r(\mathbf{x}),g(\mathbf{x}),b(\mathbf{x}))$, for $\mathbf{x} \in \mathcal{N}:=\Omega \cap \tau \mathbb{Z}^2$.
\begin{figure}[ht]
\hbox to \hsize {\hfill \epsfig{file=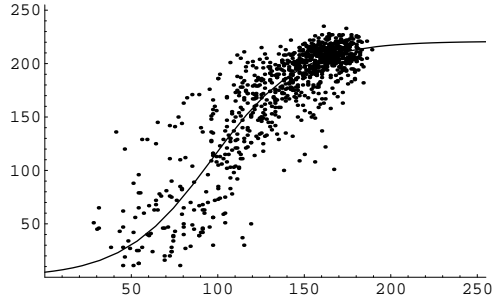,height=4cm} \hfill}
\caption{Estimation of the nonlinear curve $L$ from a distribution of points with coordinates given by the linear combination $\alpha r + \beta g +\gamma b$ of the $(r,g,b)$ color fragments (abscissa) and by the corresponding underlying gray level of the original photographs dated to 1920 (ordinate). The sensitivity parameters $\alpha,\beta, \gamma$ to the different frequencies of red, green, and blue have been chosen in order to minimize the total variance of the ordinates $\sigma^2$. Observe that, in this particular case, the curve $L$ behaves nearly linearly in the interval $[60,170]$.}
\end{figure}
Usually the gray level of an image can be described as a submanifold $\mathcal{M} \subset \mathbb{R}^3$ of dimension $1$ parameterized as $ \mathcal{M}:=\mathcal{M}_\tau=\{ \tau(gl):gl=\mathcal{L}(r,g,b):=L(\alpha r+\beta g+\gamma b),(r,g,b) \in \mathbb{R}^3_+\}$, where $\alpha,\beta,\gamma>0$, $\alpha+\beta+\gamma=1$, $L:\mathbb{R} \rightarrow \mathbb{R}$ is a non-negative increasing function, and $\tau:\mathbb{R}_+ \rightarrow \mathbb{R}^3_+$ is a section such that $\mathcal{L} \circ \tau = \text{id}_{\mathbb{R}_+}$. The function $L$ is assumed smooth, nonlinear, and normally nonconvex and nonconcave. Therefore, in our terminology, $(\mathcal{M},\mathcal{L})$ is a nonlinear projection, where $\mathcal{L}(r,g,b):=  L(\alpha r+\beta g+\gamma b)$.
For example, Figures 1 and 2 illustrate a typical situation where this model applies and Figure 4 describes the typical shape of an $L$ function, which is here estimated by fitting a distribution of data from real color fragments.

In fact, in this case, there is an area $\Omega \backslash D$ of the domain $\Omega\subset \mathbb{R}^2$ of the image, where some fragments with colors are placed and complete information samples $\mathcal{C}$ are available, and an other area $D$ where only a gray level  information $\mathcal{D}$ is known.

\begin{figure}[ht]
\hbox to \hsize {\hfill \epsfig{file=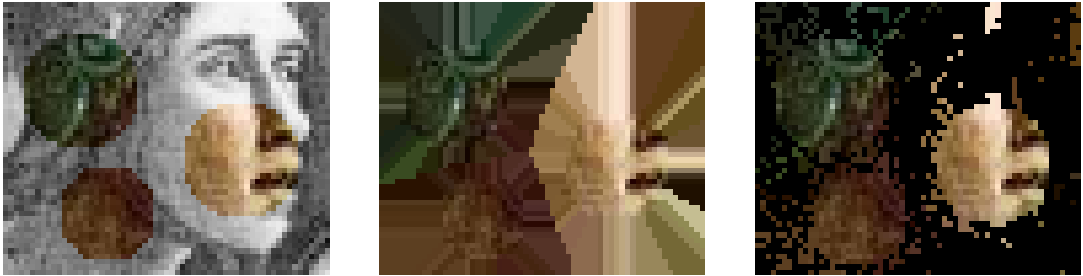,height=4cm} \hfill}
\caption{The color of the pixel is extended in the direction of the corresponding Voronoi patch. The constraint check of the projection onto the sub-manifold up to a threshold allows to keep only the color information which is compatible with the known gray information.}
\end{figure}

So our abstract Problem 1 can be reformulated as follows: When is possible to reconstruct the colors from the known information of the colors of the fragments $\mathcal{C}$ and from the gray level information $\mathcal{D}$ of the missing part?
There is not a unique answer to this question and different solutions are described depending on the methods chosen to solve the problem. Here we want to illustrate two techniques: The first based on interpolation methods and the second on PDE and variational calculus. Next we discuss how these two different strategies should be combined in order to achieve possibly better results.
\begin{figure}[hct]
\hbox to \hsize {\hfill \epsfig{file=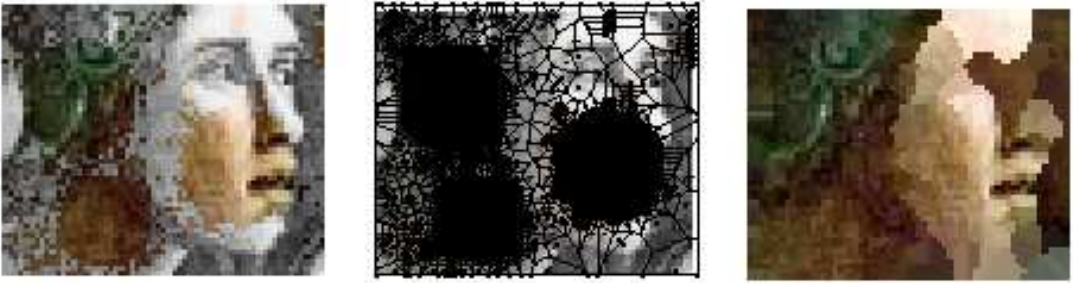,height=4cm} \hfill}
\caption{From the new deduced information, one applies again a Voronoi domain decomposition and iterates the process.}\end{figure}
\section{Nonuniform sampling methods and interpolation techniques}

On one hand if the sampling set $\mathcal{C}$ is quite spread and dense, then clearly one might apply a nonuniform sampling method, maybe those described in \cite{AG01,FS94}, or any other efficient interpolation method to recover the colors of the image. In this situation the information $\mathcal{D}$ would not turn out to be of any use. 
On the other hand, if the sampling set $\mathcal{C}$ is relatively small and concentrated, then the previous technique will fail or, at least, will give a very inaccurate result and therefore also the information $\mathcal{D}$ should be taken into account and maybe exploited.

Thus, we propose here a very simple scheme based on interpolation techniques that in the practice gives acceptable results. Assume to have at our disposal the following procedures {\bf EXTEND}, {\bf THRS}, and {\bf ESTIM}:
\\

$\mathbf{ EXTEND[F] \rightarrow G}$: given an incomplete color image {\bf F} the procedure {\bf EXTEND} computes a new color image {\bf G} which coincides with {\bf F} in the known color region of {\bf F} and it extends the colors of {\bf F} out of this region.
\\

$\mathbf{ THRS[F,F_0,\mathcal{L},\varepsilon] \rightarrow G}$: given a color image {\bf F}, a gray image $\mathbf{F_0}$, a color-gray conversion function $\mathcal{L}$, and a threshold $\varepsilon >0$ the procedure {\bf THRS} computes a new image {\bf G} which coincides with {\bf F} only in those points $\mathbf{x}$ such that $|\mathcal{L}(\mathbf{F}(\bf x))-\mathbf{F_0}(\bf x)| \leq \varepsilon$, and it is $0$ elsewhere.
\\

$\mathbf{ ESTIM[F,F_0] \rightarrow \{\mathcal L,\sigma^2\} }$: given a color image {\bf F} and a gray image $\mathbf{F_0}$ the procedure {\bf ESTIM} estimates the color-gray conversion function $\mathcal{L}$ and the total variance of the ordinates $\sigma^2$ with respect to the corresponding $L$ function as described in Figure 4.
\\

At this moment we do not care of the accuracy the {\bf EXTEND} procedure can achieve, nor its possible concrete implementation.
With these procedures at hand we define the following reconstruction algorithm:

\begin{alg}
$\mathbf{RESTORE[F,F_0]  \rightarrow G}$:\\
\noindent $\mathbf{G}=0;$\\
\noindent While $\mathbf{G} \neq \mathbf{F}$ do\\
\indent $\mathbf{G} = \mathbf{F}$;\\
\indent $\{\mathcal{L},\sigma^2\}=\mathbf{ ESTIM[F,F_0]}$;\\
\indent $\varepsilon= \Gamma(\sigma^2)$;\\
\indent $\mathbf{F=  EXTEND[F]}$;\\
\indent $\mathbf{F=  THRS[F,F_0,\mathcal{L},\varepsilon]}$;\\
\noindent od
\end{alg}

The procedure $\mathbf{RESTORE[F,F_0]  \rightarrow G}$ computes an image $\mathbf{G}$ which is the color reconstruction of the image $\mathbf{F}$ from the samples $\mathcal{C}$, being known also the gray level information $\mathcal{D}$ of the missing part $\mathbf{F_0}$. The threshold $\varepsilon>0$ used at each iteration is a function of the variance $\sigma^2$, since we cannot expect to be more precise in our reconstruction than the intrinsic uncertainty of the model describing the color-gray conversion.
\begin{figure}[ht]
\hbox to \hsize {\hfill \epsfig{file=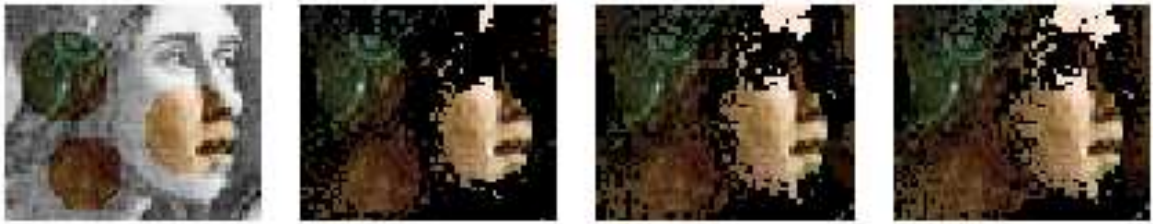,height=3.2cm} \hfill}
\caption{Successive iterations of the interpolation process. The final image is illustrated on the right. Observe that neither the nose nor the eyes of the character could have been reconstructed.}
\end{figure}
\\

In order to understand how the Algorithm 1 in fact works, let us illustrate here a simple implementation of a procedure  {\bf EXTEND}. Clearly in more concrete applications one should use more sophisticated, accurate, and robust methods \cite{AG01,FGS95,FS94}.
Given the color samples $\mathcal{C}$ of our image  $\mathbf{F}$ one splits the domain by using a Voronoi decomposition, where the nodes for the decomposition are $\mathcal{C}$, see Figure 3. One realizes a template procedure {\bf EXTEND} just extending the color of any sample in $\mathcal{C}$ in the corresponding region of the Voronoi decomposition whose it is part, see Figure 5. At this point we apply the procedure $\mathbf{THRS}$ just to keep the best possible extension which is compatible with the known gray level, up to the prescribed threshold. From this new extended set of color samples we iterate the process, constructing a new Voronoi decomposition, realizing the extension, and then applying again the thresholding and so on, see Figures 6 and 7. During the iterative scheme, a dynamical learning process of the model describing the color-gray conversion is realized.

It is not difficult to see that Algorithm 1 converges after a finite number of iterations, and that it will stop having computed some possible extension $\mathbf{G}$ of the image $\mathbf{F}$.
This $\mathbf{G}$ might be considered the best possible extension of the samples $\mathcal{C}$ compatible with the known gray level samples $\mathcal{D}$.
This method depends strongly on the way the procedure {\bf EXTEND} is implemented and it exploits the information given by $\mathcal{C}$ and $\mathcal{D}$ in an independent way. It is clear that one cannot reconstruct with this particular implementation any color which is not already included in the initial sample set $\mathcal{C}$, and not connected geometrically with one of the initial sample color pixel  by a sequence of Voronoi decompositions. Therefore, this method is highly local (i.e., no global properties of the image are in fact used) and it will fail whenever we will try to reconstruct colors which cannot depend somehow from those already known and given, and not connected geometrically with one of them. For this reason we want here to discuss also other possible inpainting techniques based on global properties of the signals.

\section{PDE and variational methods}

The modern techniques of processing digital signals are mainly and traditionally based on developing and applying Fourier and harmonic analysis concepts, e.g., wavelets \cite{Mal}, time-frequency/Gabor analysis \cite{FS,FS1,G}, sampling theory \cite{BF,BF2}, and on stochastic and Bayesian modeling \cite{WG}.
Only recently signal processing has become a very attractive field for applications of PDE and variational methods. We refer the reader to the nice introduction \cite{AK02} for a presentation of this emerging field, for major details, and an extended literature.

In this section we want first to show how Problem 1 can be modeled for 1D signals by using suitable energy functionals whose minimum corresponds to possible solutions. In particular, we illustrate a numerical scheme which can be applied in the cases where the solution exhibits enough regularity.
Next we propose a non-standard technique based on time-frequency/Gabor analysis to regularize problems involving non-smooth solutions. This method will require to move from a 1D problem to an equivalent 2D problem.
This will introduce us to the treatment of 2D digital signals (digital images) and to the formulation of a variational scheme for the solution of Problem 1 for digital images.

\begin{figure}[ht]
\hbox to \hsize {\hfill \epsfig{file=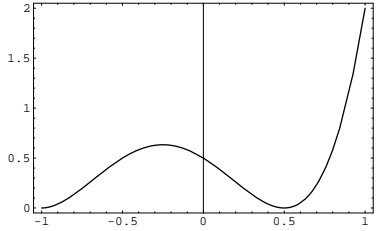,height=3cm} \hfill}
\caption{As a nonlinear distortion we consider for example $L(x)=1.8 (x+1)^2(x-\frac{1}{2})^2$.}
\end{figure}

Let us assume $\Omega \subset \mathbb{R}$, $\lambda,\mu>0$,  and $\mathcal{L} \in C^1(\mathbb{R})$. We look for solutions $u \in W^{1,p}(\Omega)$, for $p>1$, of the following
\begin{pb}
\label{pb2}
\begin{equation}
\text{arg}\inf \left \{F(u)=\mu \int_{\Omega\backslash D} |u(x)-\bar u(x)|^2 dx + \lambda \int_{D} |\mathcal{L}(u(x))-\bar u(x)|^2 dx+\int_{\Omega}|\nabla u(x)|^p dx, u \in  W^{1,p}(\Omega) \right\},
\end{equation}
where $\bar u \in W^{1,p}(\Omega)$ is the given observed signal, which is presumed to be correct on $\Omega\backslash D$ and distorted by the (nonlinear) function $\mathcal{L}$ on $D$, see Figures 8,9.
\end{pb}

\begin{figure}[ht]
\hbox to \hsize {\hfill \epsfig{file=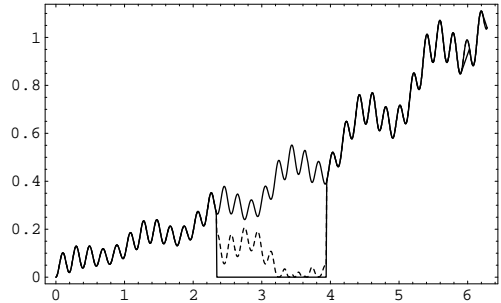,height=4cm} \hfill}
\caption{The continuous curve represents the original signal, while the dashed curve is the distorted signal by the nonlinear mapping.}
\end{figure}
To minimize the functional $F$ means essentially to find a function $u \in W^{1,p}(\Omega)$ that coincides with $\bar u$ on $\Omega \backslash D$, $\mathcal{L}(u)=\bar u$ on $D$, and it minimizes some $p$-norm of its variation. In particular, for $p=2$ we would look for solutions with small curvature. Let us fix $p=2$ in the following for simplicity.
By the {\it direct method of the calculus of variations} \cite[Section 2.1]{AK02}\cite{Ev} it is well-known that Problem 2 has a solution $u \in W^{1,2}(\Omega)$, and it solves the Euler-Lagrange equation
\begin{pb}
\label{pb3}
\begin{equation}
0=-\Delta u +2\mu(u-\bar u)1_{\Omega \backslash D} + 2 \lambda (\mathcal{L}(u)-\bar u)\frac{d\mathcal{L}}{d x}(u) 1_{D} := \mathcal{E}(\mathcal{L},u)
\end{equation}
\end{pb}
If $\mathcal{L}$ is convex then the solution of Problem 2 is unique, and Problem 3 is equivalent to Problem 2. In the case $\mathcal{L}$ is {\it not} convex then we can anyway look for solutions of Problem 2 among those of Problem 3.

\begin{figure}[ht]
\hbox to \hsize {\hfill \epsfig{file=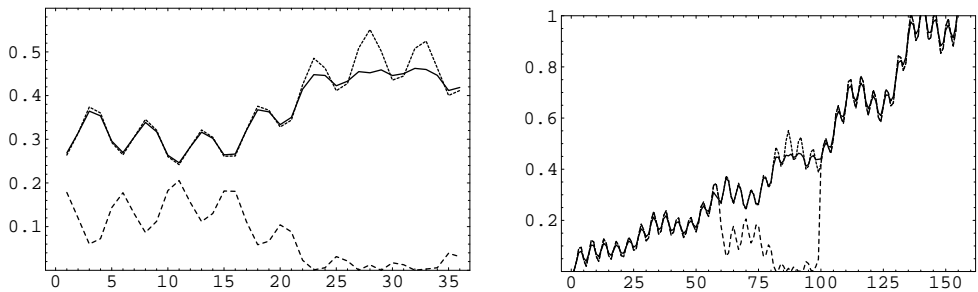,height=4cm} \hfill}
\caption{The steepest descent method is applied to recover the original signal from the partial information on the missing part and the original information on the known part. In this nice example the reconstruction of the signal is rather good.}
\end{figure}

A way to compute a solution of Problem 3 is the {\it steepest descent method}: The solution $u:\Omega \times \mathbb{R}_+ \rightarrow \mathbb{R}$ of the evolutionary equation
\begin{equation}
\label{evolution}
\frac{\partial u(x,t)}{\partial t} = -\mathcal{E}(\mathcal{L}, u(x,t)), \quad u(\cdot,0)=u_0.
\end{equation} 
tends to make vanish $\mathcal{E}(\mathcal{L}, u(x,t))$ for $t \rightarrow +\infty$. Therefore one can look for solutions  of Problem 2 given by $u(x):=\lim_{t \rightarrow \infty} u(x,t)$.

\begin{rem}
This method is as more effective as $u_0$, the initial guess, is closer to the real solution $u$.
As we will discuss later, this is the crucial and important motivation to combine interpolation methods, for achieving the best first guessing of the solution, and then variational techniques to complete the reconstruction. 
In fact the result shown in Figure 10 would have been impossible by using an interpolation (or irregular sampling) technique only, since the size of the gap (missing part) to be reconstructed would have been too large with respect to the band-width of the signal, and the use of the variational method might fail if the first guess $u_0$ is not close enough to the solution. For example, in Figure 11 we show the evolution of $u(x,t)$ in reconstructing the missing part of the signal in Figure 9, from a first linear interpolation guessing.
Therefore, we expect that in concrete and practical problems only an interlacing of these two different tools can give efficient results.
\end{rem}
\begin{figure}[ht]
\hbox to \hsize {\hfill \epsfig{file=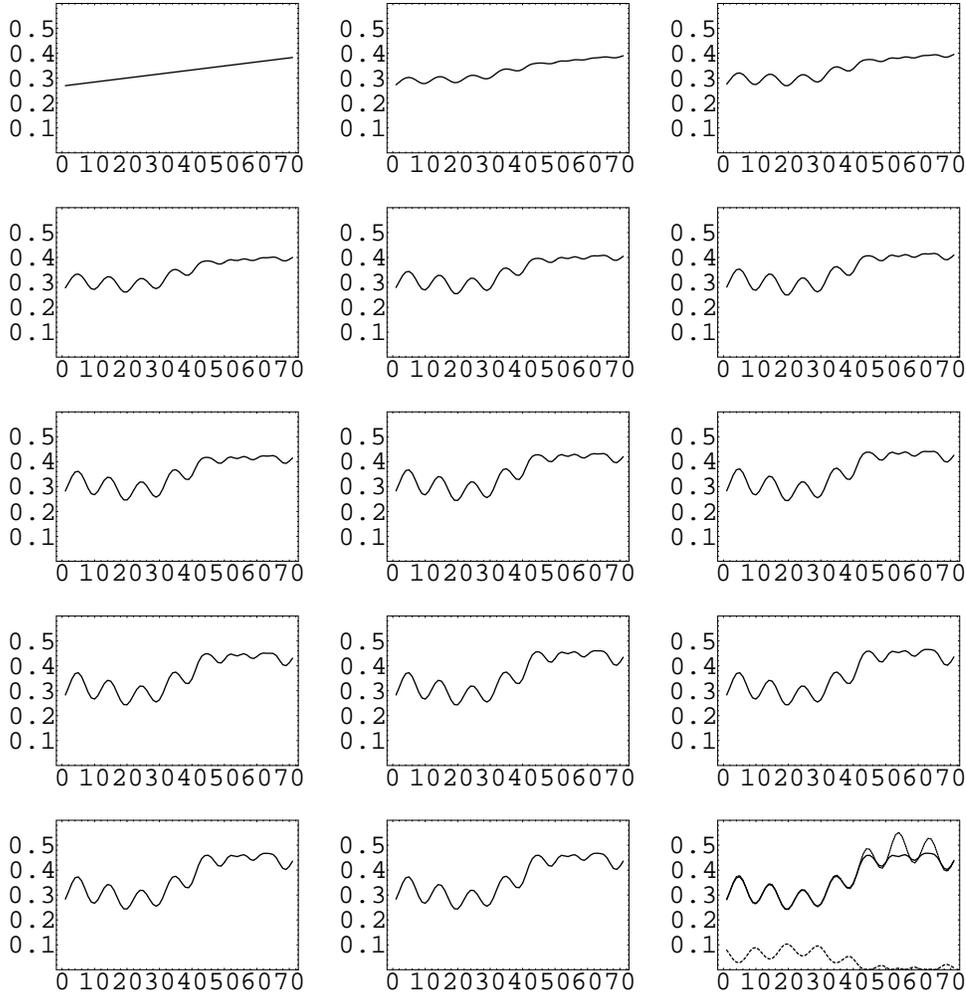,height=14cm} \hfill}
\caption{Iterations of the discretized steepest descent method}
\end{figure}
\\

The steepest descent method illustrated in \eqref{evolution} can be easily discretized and it can work on samples of the signal by means, e.g.,  of an explicit Euler scheme.
Let us denote $\Delta^2(\mathbf{v})(i):= \mathbf v(i-1)-2\mathbf v(i)+ \mathbf v(i+1)$, for any vector $\mathbf{v}$ and for $i =1,...,N-2$ and define the procedure 
$$
\mathbf{\mathcal{E}[\mathcal{L},v,i] \rightarrow \Delta^2(\mathbf{v})(i) -2\mu(\mathbf v(i) -\mathbf{\bar u}(i))\mathbf 1_{\Omega \backslash D} -2 \lambda (\mathcal{L}(\mathbf v(i)) - \mathbf{ \bar u}(i)) \frac{d \mathcal{L}}{dx}(\mathbf v(i)) \mathbf 1_{D}}, \quad i=1,...,N-2. 
$$

\begin{alg}
$\mathbf{STEEP\_DESC[u_0,\varepsilon,\mathcal{L},\Delta t]}:$\\
$\mathbf v=\mathbf{\bar v}:=\mathbf u_0$;\\
While $\text{max}\{\mathbf{\mathcal{E}[\mathcal{L},v,i]|}, i=1,...,N-2\} > \varepsilon$ do\\
\indent For i=1,...,N-2 do\\
\indent \indent $\mathbf v(i)=\overline{\mathbf v}(i)+\Delta t \left( \mathcal{E}[\mathcal{L},\overline{\mathbf v},i]\right)$;\\
\indent \indent od\\
\indent $\overline{\mathbf v} = \mathbf v$; \\
od\\
$\mathbf u = \mathbf v$.

\end{alg}
Therefore a solution of our original Problem 1 can be modeled as the result of the procedure $\mathbf{STEEP\_DESC[u_0,\varepsilon,\mathcal{L},\Delta t]}$ which, given an observed set of $N$ samples $\mathbf u_0 = \mathbf{\bar u}:=\mathcal{C} \cup \mathcal{D}$, where $\mathcal{C}$ are presumed correct and $\mathcal{D}$ are presumed distorted by the function $\mathcal{L}$, computes approximate evolutions of $ \mathbf u_0$ by discrete time steps $\Delta t$. Once the Algorithm 2 has corrected the distorted samples, one can use any interpolation result to reconstruct the solution $u$ of the Problem \ref{pb1}.
\\

A similar approach can be realized also in the case the solution is not smooth, for example $u \in BV(\Omega)$. In general the discretization method to compute a minimizing  solution $u \in BV(\Omega)$ of the Problem \ref{pb1} is more sophisticated and usually it is necessary to construct an associated regularized problem which has minimum approximating $u$.
We refer, for example, to \cite[Section 3.2]{AK02} for details on techniques used in this situation.\\

In the next section we want to present a non-standard regularization method based on time-frequency and harmonic analysis tools. 

\subsection{A non-standard regularization method: time-frequency analysis}

For a function $f \in L^2(\mathbb{R}^n)$ and $0 \neq g \in C^\infty_c(\mathbb{R}^n)$, $\|g\|_2=1$, the transformation
\begin{equation}
	V_g(f)(x,\omega):= \int_{\mathbb{R}^n} f(t) e^{-2 \pi i \omega t} \overline{g(t- x)} dt,
\end{equation}
is called the {\it short time Fourier transform} (or {\it windowed Fourier transform} or simply {\it Gabor transform}) of the function $f$.
Since $g$ is a smooth function, even if $f$ is {\it not} 
\begin{figure}[ht]
\hbox to \hsize {\hfill \epsfig{file=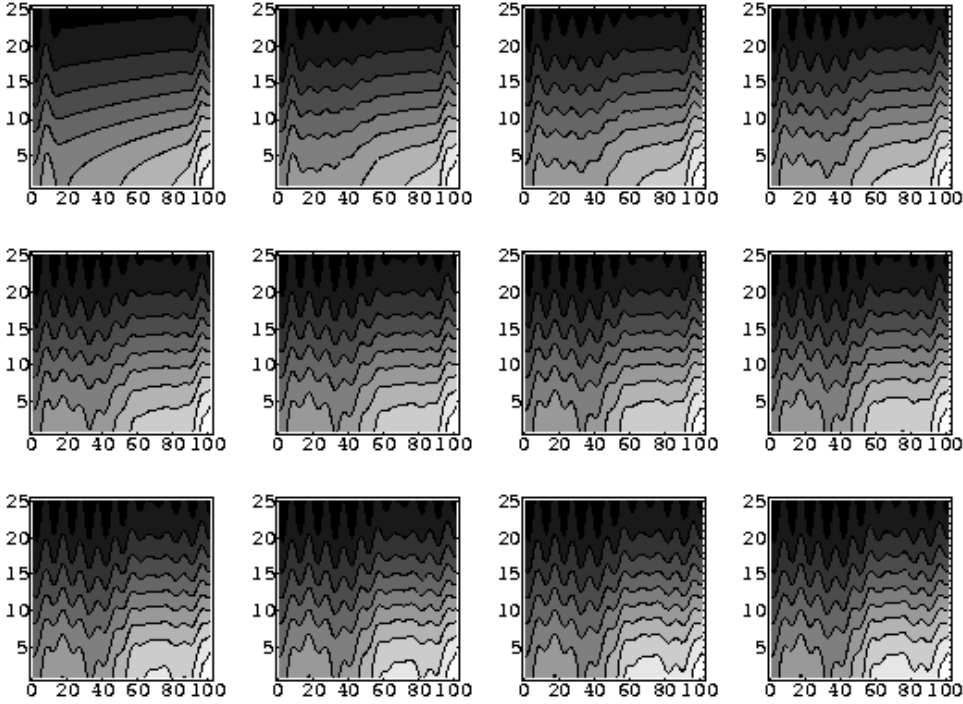,height=10cm} \hfill}
\caption{Iterations of the discretization of the steepest descent method in the time-frequency plane}
\end{figure}
smooth, $V_g(f)$ is smooth indeed. Moreover, the map $V_g: L^2(\mathbb{R}^n) \rightarrow L^2(\mathbb{R}^{2n})$ is a unitary isomorphism into its range and left-invertible by its adjoint 
\begin{equation}
V_g^*(F) = \int_{\mathbb{R}^{2n}} F(x,\omega) e^{2 \pi i \omega t} g(t- x) dx d\omega.
\end{equation}
Therefore, one has the following reproducing formula
\begin{equation}
f(t) = V_g^*\circ V_g (f)(t) =\int_{\mathbb{R}^{2n}} 	V_g(f)(x,\omega)  e^{2 \pi i \omega t} g(t- x) dx d\omega.
\end{equation}
Moreover $\Pi:=V_g \circ  V_g^*$ is the orthogonal projection from $L^2(\mathbb{R}^{2n})$ onto $\text{ran}(V_g)$.
The Gabor transform gives simultaneous information on the time/space-frequency content of a given signal. In particular, it tells us which are the ``instantaneous'' frequencies appearing in a 1D signal at a given time, somehow as the score is describing a music, telling which are the notes to be played at a given time.
\\

Because of the invertibility of $V_g$, to analyze or synthetize a signal in the time, or to do that in the time-frequency plane are equivalent operations.
For example, in Figure 12 we have visualized the Gabor transform of the iterations in Figure 11.
This suggests that Problem \ref{pb2} can be reformulated in the time-frequency plane.
In fact, as already mentioned, it is known that $u \in L^2(\Omega)$ implies $V_g(u) \in L^2(\Omega+\supp(g)\times \mathbb{R}) \cap C^\infty$.
The treatment of the problem as defined in the time-frequency plane is non-standard, especially because the domain itself is unbounded.
On the other hand, the range of $V_g$ is always constituted by smooth functions and this of course is a way to regularize the problem.
\begin{pb}
\label{pb5}
$$
\text{arg}\inf  \{F(U)=\mu \int_{\Omega_g\backslash D_g \times \mathbb{R}} |U(x,\omega)-V_g(\bar u)(x,\omega)|^2 dx d\omega + \lambda \int_{D \times \mathbb{R}} |\mathcal{L}_{TF}(U)(x,\omega)-V_g(\bar u)(x,\omega )|^2 dxd\omega
$$
\begin{equation}
+\int_{\Omega_g\times \mathbb{R}}| \nabla U|^2 dxd\omega, U \in  W^{1,2}(\Omega_g \times \mathbb{R})\},
\end{equation}
where $\Omega_g = \Omega + \supp(g)$, $D_g=D+\supp(g)$, $\mathcal{L}_{TF}(U):= V_g(\mathcal{L}(V_g^*(U)))$, and $\bar u \in L^2(\Omega)$ is the given observed signal, which is presumed to be correct on $\Omega\backslash D$ and distorted by the (nonlinear) function $\mathcal{L}$ on $D$.
\end{pb}
As usual, one can correctly compute the associated Euler-Lagrange equations and to formulate the corresponding steepest descent method, in the same way as we have done in the 1D case.
Since similar definitions and properties of the Gabor transform can be formulated for discrete signals one can easily discretize the steepest descent method and define a suitable numerical scheme. We refer the reader to \cite{FS,FS1} for major details on numerical Gabor analysis in the treatment of digital signals. The analysis of Problem \ref{pb5} and its connections with corresponding properties in the time domain will be investigated in successive contributions.
\\

This section has been useful to us also to move our attention from the 1D situation back to a 2D problem, which was the inspiration of this paper. In the next section we want to generalize what we have discussed for univariate signals on the real line to multivariate signals on 2D domains. 
In particular, it is well known that the Laplacian appearing in the Euler-Lagrange equations associated to Problem \ref{pb5} has very strong isotropic smoothing properties and does not preserve edges.
While this could not be a relevant problem for 2D function in the range of $V_g$, since they are smooth and affected by the Heisenberg uncertainty principle which makes them intrinsically ``blurred'', it is not suitable for the reconstruction of 2D functions representing natural images, which are usually characterized by the presence of edges, curves, and maybe fractal structures.

\subsection{Variational models for image inpainting}

There are several and different variational methods for solving the so called {\it image inpainting problem}, i.e., the reconstruction of a small missing part of an image by using the information of the remaining relevant known part. Each of them offers some nice properties and effects together with some drawbacks. In this section we want to propose a modification of what we consider the most simple and well known approach, in order to have a reasonable numerical solver for our original problem of the reconstruction of colors.
First, let us describe in the following three of the models discussed recently in the literature:
\begin{itemize}
\item[1.] The Beltramio-Shapiro-Caselles-Ballester model. Let $L(u)$ be a smoothness measure of the image $u$, for example, a second order differential operator given by
\begin{equation}
L(u):=f(\nabla u, \nabla \otimes \nabla u).
\end{equation}

Typical choice for $L$ is the Laplacian $L(u)=\Delta u=\text{trace}(\nabla \otimes \nabla u)$.
The reconstruction model proposed by Beltramio, Shapiro {\it et al.} \cite{BSCB01} is based on the evolutionary equation
\begin{equation}
\frac{\partial u((x,y),t)}{\partial t}=\nabla^{\perp} u((x,y),t) \cdot \nabla L(u(x,y),t)),
\end{equation}
where $\nabla^{\perp} u=(-\frac{\partial u((x,y),t)}{\partial y},\frac{\partial u((x,y),t)}{\partial x})=|\nabla u| \mathbf{T}$, point to the tangent $\mathbf{T}$. This means that, at the equilibrium, one has $\frac{\partial  L(u)}{\partial \mathbf{T}} =0$. Thus, in terms of boundary smoothness data, the inpainting process evolves transporting the information along the extended {\it isophotes} (i.e., level curves)into the region to be reconstructed. This method has a nice numerical implementation in \cite{BSCB01} which allowed the authors to show several examples where it works nicely. 
However, since there is not a clear information connecting different isophotes, the evolution might create nonexistent T-junctions inside the inpainting region and smoothing operations are necessary to connect properly such discontinuities. Moreover, a rigorous theoretical description of such equation is still a matter of investigation.
\item[2.] The Chan and Shen non-linear diffusion model. This model is much more in the spirit of the energy methods presented in the 1D model and it offers several advantages. In particular, the mathematics related to such model is much more developed and numerical implementations have been recently well formulated. It is based on the following minimizing problem \cite{CS02}
\begin{equation}
\text{arg}\inf\left\{F(u)=\int_{\Omega\backslash D} |u(x)-\bar u(x)|^2 dx + 
\int_{\Omega}\phi(|\nabla u(x)|) dx, u \in  W^{1,p}(\Omega) \right\},
\end{equation}
where $\bar u \in BV(\Omega)$ is the given observed signal, which is presumed to be correct on $\Omega\backslash D$ and $D$ is the region to be reconstructed.

It is known that a solution of the problem (see for example \cite[Proposition 3.3.4, Section 3.3.1]{AK02}) that a minimizing solution solves in particular the following Euler-Lagrange equation
\begin{equation}
\label{cs}
0=-\nabla \cdot \left(\frac{\phi'(|\nabla u|)}{|\nabla u|}\nabla u\right)+ 2(u(x)-\bar u(x))1_{\Omega \backslash D},
\end{equation} 
endowed with suitable boundary conditions, where $\nabla u$ is the {\it approximate derivative} of $u$.

On one hand, such an evolutionary equation and the corresponding steepest descent method are normally used for the numerical solution of the energy problem as reasonable practical methods. On the other hand, this procedure is not completely mathematically correct and more sophisticated techniques should be implemented and we refer the reader to \cite[Section 3.2.3-4]{AK02} for more details. Such techniques are based again on regularization methods.
It has been also proposed by Chan and Shen \cite{CS02} to modify equation \eqref{cs} in the following way
\begin{equation}
\label{mod}
0=-|\nabla u| \nabla \cdot \left(\frac{\phi'(|\nabla u|)}{|\nabla u|}\nabla u\right)+ 2(u(x)-\bar u(x))1_{\Omega \backslash D}.
\end{equation} 
The steepest descent method in the inpainting region $D$, for the trivial choice $\phi(t)=t$, becomes
\begin{equation}
\frac{\partial u}{\partial t} = |\nabla u| \nabla \cdot \left(\frac{\nabla u}{|\nabla u|}\right),
\end{equation}
which describes the well-known (morphological invariant) {\it mean curvature motion}, see, e.g.,  \cite{ALM,MCM,CS02}. Efficient numerical implementations of such an evolutionary equation have been recently formulated for example in \cite{KWBE}. Such model tends to approximate the isophotes in the inpainting region by straight lines and, in the following, we will show that one can modify it to avoid this problem. 
\item[3. ] The Ambrosio and Masnou model. A method to extend in a smooth way the isophotes from T-junctions at the border of the inpainting region has been recently proposed and studied by  Ambrosio and Masnou \cite{AM03} based on the following minimization problem
\begin{equation}
\text{arg}\inf \left \{F(u)=\int_{\tilde D}  | \nabla u| \left (1+ \left | \frac{ \nabla u}{| \nabla u|}\right |^p \right) dx, u \in  W^{1,p}(\Omega) \right\},
\end{equation}
where $\tilde D$ is a slightly larger domain containing the inpainting region $D$. 
In fact, in their paper \cite{AM03} they discuss the corresponding relaxed functional $\bar F$ with respect to the $L^1$ convergence and the existence of $BV(\Omega)$ minimizers of $\bar F$.
In the case $n=2$ and $p>1$ Chan and Shen derived the (fourth order) Euler-Lagrange equation corresponding to such energy problem \cite{CS02}. 
\end{itemize}
From what we have just illustrated above, there are several nice intuitions and models based on variational methods in order to solve the inpainting problem, where deep mathematics and numerical problems are involved and still rather unexplored. In our understanding, presently the mathematically most well supported and the most promising in terms of efficient numerical schemes and realizations is the Chan-Shen approach.
We want now to borrow their approach in order to discuss and to model our color reconstruction problem.
We assume then that $u$ is a RGB image.
\begin{pb}
\label{pb6}
\begin{equation}
\text{arg}\inf \left \{F(u)=\mu\int_{\Omega\backslash D} |u(x)-\bar u(x)|^2 dx + \lambda \int_{D} |\mathcal{L}(u(x))-\bar u(x)|^2 dx+\int_{\Omega}\phi(|\nabla u(x)|) dx, u \in  W^{1,p}(\Omega,\mathbb{R}^3_+) \right\},
\end{equation}
where $\bar u$ is the given observed image, which is presumed to be colored on $\Omega \backslash D$ and just gray level on $D$ where the transformation from RGB to gray level is given by the (nonlinear) function $\mathcal{L}$. 
\end{pb}
Since we are dealing with natural images one should relax the functional and look for BV minimizers. This strategy can be handled in a similar way as it is done, for example, in \cite[Theorems 3.2.1-2]{AK02} and it will be discussed elsewhere. However it is true that the BV solution of the relaxed formulation of Problem \ref{pb6} is anyway solving the following Euler-Lagrange equations

\begin{equation}
\label{stat}
0=-\nabla \cdot \left(\frac{\phi'(|\nabla u_k|)}{|\nabla u_k|}\nabla u_k\right)+2\mu(u_k-\bar u_k)1_{\Omega \backslash D} + 2 \lambda (\mathcal{L}(u)-\bar u)\frac{\partial \mathcal{L}}{\partial x_k}(u) 1_{D} := \mathcal{E}_k(\mathcal{L},u), \quad k=1,2,3,
\end{equation}
where $u_1:=r, u_2:=g, u_3:=b$ are the RGB components of the image $u$, and $\nabla u_k$ is the approximate derivative.
One can consequently formulate the three evolutionary steepest descent equations by
\begin{equation}
\label{evolution2}
\frac{\partial u_k((x,y),t)}{\partial t} = -\mathcal{E}_k(\mathcal{L}, u((x,y),t)), \quad u(\cdot,0)=u_0, k=1,2,3,
\end{equation} 
to make vanish $\mathcal{E}_k(\mathcal{L}, u((x,y),t))$, $k=1,2,3$, for $t \rightarrow +\infty$. Therefore one can look for solutions  of Problem \ref{pb6} given by $u(x,y):=\lim_{t \rightarrow \infty} u((x,y),t)$.
Observe that, in particular, at the equilibrium and on the inpainting region $D$ one has
$$
0=-\nabla \cdot \left(\frac{\phi'(|\nabla u_k|)}{|\nabla u_k|}\nabla u_k\right)+ 2 \lambda (\mathcal{L}(u)-\bar u)\frac{\partial \mathcal{L}}{\partial x_k}(u).
$$
This means that the extension of the isophotes will not be in general just straight lines, as we have already seen also in the 1D model (Figure 10) where the reconstructed part is not just a linear interpolation as one could expect using a mean curvature motion evolutionary equation. 
\\

It is important to observe that \eqref{stat} and \eqref{evolution2} are systems of coupled second order Hamilton-Jacobi equations and the analysis of the solutions constitutes itself an open problem of independent interest that will be discussed elsewhere. Some results on existence and uniqueness of viscosity solutions for systems of steady nonlinear second order PDE have been proposed by Koike \cite{K1,K2,K3}. Unfortunately those results apply under \emph{quasi-monotonicity conditions} that here could not in general be verified. However, one can easily construct solutions of \eqref{stat} and \eqref{evolution2} in particular cases as follows.
\\

Let us assume that a gray level image $v$ can be identified with an RGB image $\tau \circ v$, by means of a section $\tau:\mathbb{R}_+ \rightarrow \mathbb{R}^3_+$ such that $\mathcal{L} \circ \tau = \text{id}$. For example, if $\mathcal{L}(r,g,b)=\alpha r+\beta g+\gamma b$, with $\alpha+\beta+\gamma=1$ then one can choose $\tau(gl)=(gl,gl,gl)$. 
Suppose that $u$ is an image with piecewise linear/straight isophotes, colored on $\Omega \backslash D$ and gray on $D$. Therefore we can always assume that $u=(u_1,u_2,u_3)=u_{|\Omega \backslash D} +\tau(u_{|D})$ is a vector valued function. For the choice $\phi(t)=t$, it is not difficult to show that such a function $u$ is almost everywhere a stationary solution of \eqref{evolution2}, i.e., for the choice $u_0 = \bar u:=u$ there will be no evolution, and $u$ is solution of \eqref{stat}.
\begin{figure}[hct]
\hbox to \hsize {\hfill \epsfig{file=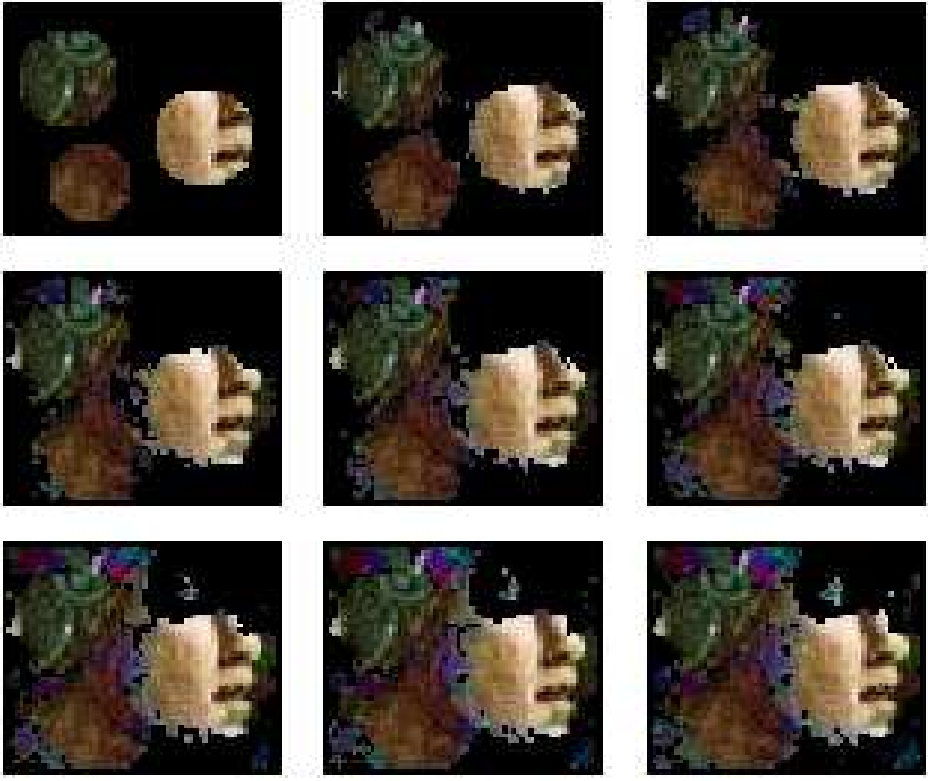,height=11cm} \hfill}
\caption{Successive iterations of the nonlinear diffusion process. Starting from Figure 2 the color is diffused extending the isophotes in the inpainting region, following the tracks of the underlying gray level. Emblematic is the reconstruction of the color on the nose of the character, where the diffusion progresses along it, but not through its boundary, the isophotes having there almost vanishing curvature.
The last image on the bottom-right is the result of 300 iterations of the outer loop of Algorithm 3. Here we have chosen $\Delta t = 0.1$, $\lambda=\mu=10$, and $\mathcal{L}(r,g,b) = \frac{1}{3}(r+b+g)$. The initial information content does not appear sufficient for a correct color reconstruction and only a limited part of the image can be partially restored.}
\end{figure}

\subsection{The numerical implementation} 

The steepest descent method illustrated in \eqref{evolution2} can be easily discretized and it can work on samples of the image by means, e.g.,  of an explicit Euler scheme. We assume here $\phi(t)=t$ and we modify again the velocity of the evolution as in \eqref{mod} by multiplication of the first term with the gradient. This produces a motion of mean curvature of the isophotes of the images.\\

For a color digital image $\mathbf{v}=(\mathbf{v}_1,\mathbf{v}_2,\mathbf{v}_3)$ and for $i =1,...,N-2$, $j=1,...,M-2$, let us denote
\begin{equation}
\kappa(\mathbf{v}_k)(i,j):= \sum_{(m,n) \in \mathcal{N}(i,j)} \frac{2 ( \mathbf{v}_k(m,n)- \mathbf{v}_k(i,j))}{|\nabla \mathbf{v}_k(i,j)| +|\nabla \mathbf{v}_k(m,n)|} 
\end{equation}
where $\mathcal{N}(i,j)$ denoted the 4-neighborhood of the pixel position $(i,j)$, and $|\nabla \mathbf{u}(i,j)| = \sqrt{(\delta_{x}^c \mathbf{u}(i,j))^2 + (\delta_{y}^c \mathbf{u}(i,j))^2}$, $\delta_x^c, \delta_y^c$ being the central differences of order 1 in the direction $x$ and $y$ respectively. Then define the procedure 
\begin{eqnarray*}
\mathbf{\mathcal{E}_k[\mathcal{L},v,i,j]} &\rightarrow& \mathbf{|\nabla \mathbf{v}_k(i,j)|\kappa(\mathbf{v}_k)(i,j) -2\mu(\mathbf v_k(i,j) -\mathbf{ \bar u_{k}}(i,j))\mathbf 1_{\Omega \backslash D} -2 \lambda (\mathcal{L}(\mathbf v(i,j)) - \mathbf{ \bar u}(i,j)) \frac{d \mathcal{L}}{dx_k}(\mathbf v(i,j)) \mathbf 1_{D}}, \\
&& i=1,...,N-2, \quad j=1,...,M-2, \quad k=1,2,3.
\end{eqnarray*}

\begin{alg}
$\mathbf{STEEP\_DESC\_2D[u_0,\varepsilon,\mathcal{L},\Delta t]}:$\\
$\mathbf v:=\overline{\mathbf v}:=\mathbf u_0$;\\
While $\text{max}\{|\mathbf{\mathcal{E}[\mathcal{L},v,i,j]|}, i=1,...,N-2, j=1,...,M-2\} > \varepsilon$ do\\
\indent For i=1,...,N-2 do \\
\indent \indent  For j=1,...,M-2 do \\
\indent \indent  \indent For k=1,...,3 do \\
\indent \indent \indent \indent $\mathbf v_k(i,j)=\overline{\mathbf v}_k(i,j)+\Delta t \left( \mathcal{E}_k[\mathcal{L},\overline{\mathbf v},i,j]\right)$;\\
\indent \indent \indent od \\
\indent \indent \indent $\overline{\mathbf v}= \mathbf v$; \\
\indent \indent od\\
\indent od \\
od\\
$\mathbf u = \mathbf v$.

\end{alg}

\begin{figure}[H]
\hbox to \hsize {\hfill \epsfig{file=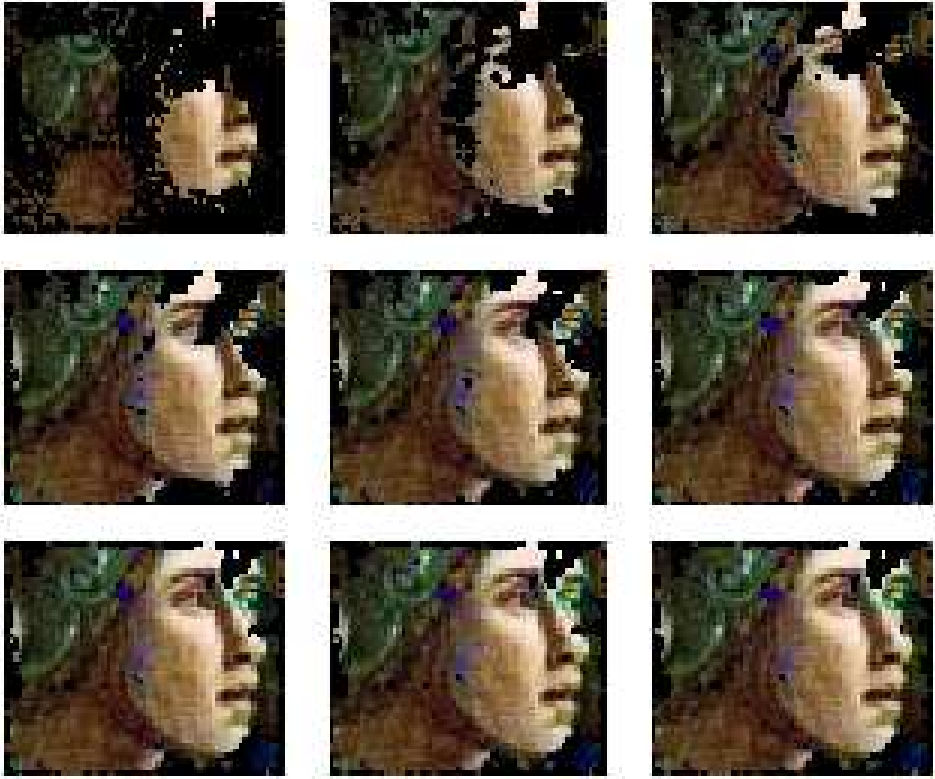,height=11cm} \hfill}
\caption{Successive iterations of the nonlinear diffusion process. Starting from a few colored samples more the diffusion process is able almost to complete the reconstruction of the color of all the face of the character after 300 iterations. Again we have chosen here $\Delta t = 0.1$, $\lambda=\mu=10$, and $\mathcal{L}(r,g,b) = \frac{1}{3}(r+b+g)$. }
\end{figure}

Therefore a solution of our original Problem 1 can be modeled as the result of the procedure $\mathbf{STEEP\_DESC\_2D[u_0,\varepsilon,\mathcal{L},\Delta t]}$ which, given an observed set of $N \times M$ samples $\mathbf u_0 = \mathbf{\bar u}:=\mathcal{C} \cup \mathcal{D}$, where $\mathcal{C}$ are presumed colored and $\mathcal{D}$ are presumed gray, computes approximate evolutions of $ \mathbf u_0$ by discrete time steps $\Delta t$.
Of course the explicit Euler scheme is not the most efficient, and implicit variants can be easily derived, see, e.g., \cite{KWBE}.

Let us show now some applications of Algorithm 3 for the color inpainting problem in Figure 13 and Figure 14.


\subsection{Conclusion: The interpolation-inpainting method}

Our reconstruction method is based on an inpainting model which exploits additional information given by a known nonlinear projection of the solution onto a manifold of lower dimension. The numerical solution can be implemented by a suitable discretization of evolutionary equations which are assumed to start from an initial guess of the solution.
We have shown that in concrete problems where the inpainting region is large with respect to the known complete data, only a combination of suitable interpolation and variational methods can work properly. Since the interpolation methods can be computationally expensive (for example the construction of a Voronoi decomposition is indeed rather expensive) one can combine a few initial iterations of Algorithm 1, in order to extend as much as possible the knowledge about the color-gray conversion and a possibly good initial guess of the solution, with successive fast iterations of Algorithm 3 to complete the reconstruction, see Figure 15.

\begin{figure}[H]
\hbox to \hsize {\hfill \epsfig{file=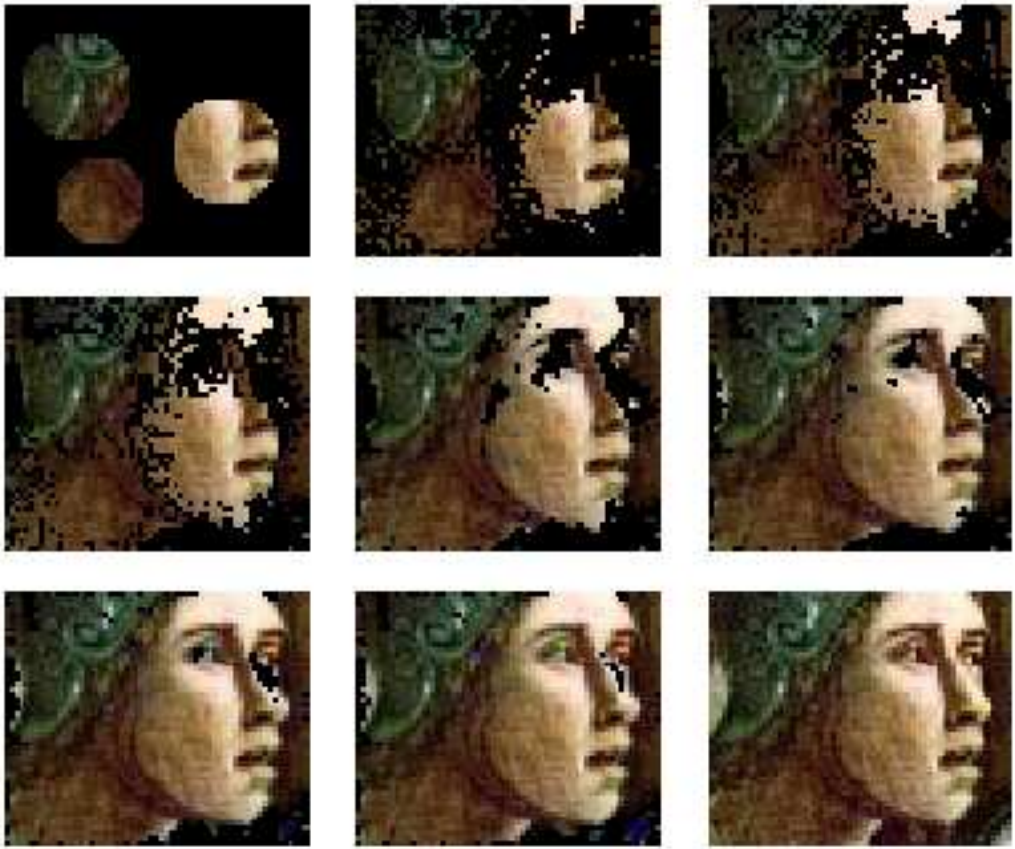,height=12cm} \hfill}
\caption{Successive iterations of the interpolation-inpainting process. The first three pictures on the top are the result of iterations of Algorithm 1, and the successive are generated by Algorithm 3. Again we have chosen here $\Delta t = 0.1$, $\lambda=\mu=10$, and $\mathcal{L}(r,g,b) = \frac{1}{3}(r+b+g)$. }
\end{figure}

\bibliography{mfornasiref}

\bigskip

\noindent Massimo Fornasier, Universit\`a ``La Sapienza'' in Roma, 
Dipartimento di Metodi e Modelli Matematici per le Scienze Applicate\\
Via Antonio Scarpa, 16/B, I-00161 Roma, Italy.\\
and\\
Universit\"at Wien, Institute F\"ur Mathematik, NuHAG \\
Nordbergstrasse 15, A-1090 Wien, Austria.\\
email: {\tt mfornasi@math.unipd.it}\\

\newpage
\section*{Vitae}

{\bf Massimo Fornasier} received his Ph.D degree in Computational Mathematics on February 2003 at the University of Padova, Italy. Within the European network RTN HASSIP (Harmonic Analysis and Statistics for Signal and Image Processing) HPRN-CT-2002-00285, he cooperated as PostDoc with NuHAG (the Numerical Harmonic Analysis Group), Faculty of Mathematics of the University of Vienna, Austria and the AG Numerical/Wavelet-Analysis Group of the Department of Mathematics and Computer Science of the Philipps-University in Marburg, Germany (2003). Since June 2003 he is research assistant at the Department of Mathematical Methods and Models for the Applied Science at the University ``La Sapienza'' in Rome. Since May 2004 he is Individual Marie Curie Fellow (project FTFDORF-FP6-501018) at NuHAG.
His research interests include applied harmonic analysis with particular emphasis on time-frequency analysis and decompositions for applications in signal and image processing. Since 1998, he developed with Domenico Toniolo the Mantegna Project (http://www.pd.infn.it/$\sim$labmante/) at the University of Padova and the local laboratory for image processing and applications in art restoration.
Recently he has focused his attention on adaptive and dynamical schemes for the numerical solution of (pseudo)differential equations and inverse problems in digital signal processing.
\end{document}